\documentclass[11pt,a4paper]{article}
\usepackage{amssymb, amsmath, latexsym, proof}

\usepackage{xcolor,soul}

\newtheorem{theorem}{Theorem}
\newtheorem{lemma}{Lemma}[section]

\newtheorem{corollary}[lemma]{Corollary}
\newtheorem{proposition}[lemma]{Proposition}

\newtheorem{definition}{Definition}
\newtheorem{example}[lemma]{Example}
\newtheorem{remark}[lemma]{Remark}

\newcommand{\bl}{\begin{lemma}}
\newcommand{\el}{\end{lemma}}
\newcommand{\bt}{\begin{theorem}}
\newcommand{\et}{\end{theorem}}
\newcommand{\bcor}{\begin{corollary}}
\newcommand{\ecor}{\end{corollary}}
\newcommand{\bp}{\proof{.}}
\newcommand{\ep}{\eop}
\newcommand{\bpr}{\begin{proposition}}
\newcommand{\epr}{\end{proposition}}
\newcommand{\brem}{\begin{remark} \em}
\newcommand{\erem}{\end{remark}}
\newcommand{\bd}{\begin{definition} \em}
\newcommand{\ed}{\end{definition}}
\newcommand{\bex}{\begin{example} \em
}
\newcommand{\eex}{\end{example}}
\newcommand{\beq}{\begin{equation} }
\newcommand{\eeq}{\end{equation}}

\newcommand{\bi}{\begin{itemize}
  }
\newcommand{\ei}{\end{itemize}}
\newcommand{\ben}{\begin{enumerate} }
\newcommand{\een}{\end{enumerate} }

\newcommand{\refeq}[1]{(\ref{#1})}
\newenvironment{enumr}{

\begin{enumerate}     }{\end{enumerate}

}

\newcommand{\benr}{\begin{enumr}
  }
\newcommand{\eenr}{
\end{enumr}}

\newcommand{\ignore}[1]{}

\newcommand{\al}[1]{\forall #1\:}
\newcommand{\ex}[1]{\exists #1\:}

\newlength{\hilflh}

\newcommand{\naturals}{\mathbb{N}}

\renewcommand{\emptyset}{\varnothing}

\newcommand{\cL}{{\mathcal L}}

\newcommand{\cF}{{\mathcal F}}

\newcommand{\fG}{\mathfrak{G}}

\newcommand{\ga}{\alpha}
\newcommand{\gb}{\beta}
\newcommand{\gd}{\delta}
\renewcommand{\ge}{\varepsilon}
\newcommand{\gl}{\lambda}
\newcommand{\gL}{\Lambda}

\newcommand{\gy}{\gamma}
\newcommand{\gw}{\omega}
\newcommand{\gS}{\Sigma}

\renewcommand{\phi}{\varphi}

\newcommand{\imp}{\rightarrow}
\newcommand{\eqv}{\leftrightarrow}

\newcommand{\Tr}{\mathrm {Tr}}
\newcommand{\LL}{\mbox{\textit{\textbf{PL}}}}

\newcommand{\Glp}{\mathrm{GLP}}

\newcommand{\Con}{\mathrm{Con}}

\newcommand{\Prf}{\mathrm{Prf}}

\newcommand{\PRA}{\mathsf{PRA}}

\newcommand{\RFN}{\mathsf{RFN}}

\newcommand{\PA}{\mathsf{PA}}
\newcommand{\EA}{\mathrm{EA}}

\newcommand{\ACA}{\mathsf{ACA}}

\newcommand{\rst}{\upharpoonright}
\newcommand{\gn}[1]{\ulcorner #1 \urcorner}

\newcommand{\ord}[2]{|#2|_{\Pi_1^0}}

\newcommand{\la}{\langle}
\newcommand{\ra}{\rangle}

\renewcommand{\models}{\vDash}      

\newcommand{\On}{\mathrm{On}}
\newcommand{\Cr}{\mathrm{Cr}}
\newcommand{\en}{\mathrm{en}}

\newcommand{\Lim}{\mathrm{Lim}}

\newcommand{\iffdef}{\stackrel{\text{def}}{\iff}}

\newcommand{\nat}{\naturals}

\renewcommand{\leq}{\leqslant}
\renewcommand{\geq}{\geqslant}

\newcommand{\Rc}{\mathrm{RC}}

\newcommand{\dom}{\text{dom}}

\newcommand{\btR}{\overline{\mathsf{R}}}

\newcommand{\eop}{$\Box$ \protect\par \addvspace{\topsep}}
\newcommand{\proof}[1]{\protect\par\addvspace{\topsep}\noindent {\bf Proof#1}}

\newcommand{\UTB}{\mathsf{UTB}}
\newcommand{\CT}{\mathsf{CT}}
\newcommand{\TT}{\mathsf{T}}
\newcommand{\num}{\underline}
\newcommand{\Set}{\mathrm{Set}}

\newcommand{\tRFN}[2]{{#1}\textrm{-}\RFN({#2})}
\newcommand{\tCA}{\textsf{-CA}}
\newcommand{\lh}{\mathrm{lh}}
\newcommand{\bit}{\mathrm{bit}}

\renewcommand{\LL}{\mathcal{L}}
\newcommand{\tR}{\mathsf{R}}

\newcommand{\bgvee}{\textstyle\bigvee}

\newcommand{\Wo}{\mathbb{W}}
\newcommand{\Fo}{\mathbb{F}}

\newcommand{\Rcl}{{\Rc_\gL}}

\newcommand{\IB}{\mathsf{IB}}
\newcommand{\ATR}{\mathsf{ATR}}

\renewcommand{\ord}{\mathrm{ord}}

\begin{document}

\title{Conservativity spectra and generalized Ignatiev model\thanks{This work was performed at the Steklov International Mathematical Center and supported by the Ministry of Science and Higher Education of the Russian Federation (agreement no. 075-15-2019-1614).}}
\author{Lev D. Beklemishev\footnote{The corresponding author,  e-mail: \texttt{bekl@mi-ras.ru}.} \\
Steklov Mathematical Institute of Russian Academy of Sciences\\ Gubkina str. 8, 119991 Moscow, Russia}
\maketitle
\begin{abstract}
We study a generalization of the notion of conservativity spectrum of an arithmetical theory to a language with transfinitely many truth definitions introduced in \cite{BekPakh}. We establish a correspondence of conservativity spectra and points of a generalized Ignatied model introduced and studied by D. Fern\'andez-Duque and J. Joosten. We also show that the results of \cite{BekPakh} easily yield the so-called Schmerl formulas for iterated reflection principles of predicative strength.
\end{abstract}

\section{Introduction}

This paper is a postscript to \cite{BekPakh} where methods of provability logic and reflection algebras were applied to the proof-theoretic analysis of theories of predicative strength. In particular, the notion of \emph{conservativity spectrum} of an arithmetical theory was generalized to the languages in which sets of the hyperarithmetical hierarchy and the corresponding truth definitions are expressible. Very roughly, the conservativity spectrum of a theory $S$ is a transfinite sequence of ordinals $\gb=\ord_\ga(S)$ such that the $\gb$-th iteration of the reflection principle for the formulas of complexity $\Pi_{1+\ga}$ over some fixed base theory is contained in $S$. Thus, the conservativity spectrum of a theory $S$ carries the information about the strength of $S$ with respect to formulas in each complexity class $\Pi_{1+\ga}$.

Conservativity spectra corresponding to the levels of the arithmetical hierarchy\footnote{More precisely, what we call here $\gw$-conservativity spectra.} were introduced by J.~Joosten \cite{Joo15a} under the name \emph{Turing--Taylor expansions} of arithmetical theories. He established a one-to-one correspondence between $\gw$-conservativity spectra and points of the so-called Ignatiev Kripke frame. Later,
it was shown in~\cite{Bek18} that the Ignatiev frame can also be seen as a natural algebraic model of the reflection calculus extended by the conservativity operators.
Ignatiev's frame was originally introduced in \cite{Ign93} as a universal model for the variable-free fragment of Japaridze's provability logic GLP. 

In this paper we extend the results of \cite{Joo15a} to the more general notion of $\lambda$-conservativity spectrum, for any constructive ordinal $\gl$, as defined in \cite{BekPakh}. Fern\'andez and Joosten \cite{FJ13a} introduced an extension of the Ignatiev frame to the language with transfinitely many modalities. Our main theorem shows that the worlds of this frame, the so-called $\ell$-sequences in the terminology of \cite{FJ13a}, coincide with $\gl$-conservativity spectra. By \cite{FJ13a}, $\ell$-sequences have a simple characterization in terms of ordinal functions related to the Veblen hierarchy. Thus, our result gives an explicit answer to the basic question what sequences of ordinals can actually occur as conservativity spectra of theories and confirms a natural conjecture of \cite{BekPakh}.

\ignore{
One of the basic questions about conservativity spectra is what sequences of ordinals can actually occur as conservativity spectra of theories? In this paper we settle a conjecture of \cite{BekPakh} characterizing such sequences for the class of extensions of elementary arithmetic. As conjectured, these sequences coincide with the so-called $\ell$-sequences of ordinals studied by J. Joosten and D. Fern\'andez \cite{FJ13a}. The $\ell$-sequences have a simple characterization in terms of ordinal functions related to the Veblen hierarchy. They constitute the set of worlds in the so-called extended Ignatiev Kripke frame. Ignatiev's frame was introduced in \cite{Ign93} as a universal frame for the variable-free fragment of Japaridze's provability logic GLP. It has been somewhat modified and studied in \cite{BJV,Ica09}. Fern\'andez and Joosten \cite{FJ13a} introduced an extension of the Ignatiev frame to the language with transfinitely many modalities.
}

This paper is not self-contained and borrows a lot of material, notations and results from \cite{BekPakh}. Therefore, we presuppose the reader's familiarity with that paper. However, the results coming from other sources are explained here in more detail.

\section{Preliminaries}
Here we briefly summarize the framework of iterated truth theories to which our results apply. We follow the presentation in \cite{BekPakh} where additional details can be found.

\subsection{Iterated Tarskian truth definitions}
Let $\LL_0$ denote the language of elementary arithmetic $\EA$. We fix an elementary definable well-ordering $(\Lambda,<)$ and consider the language $$\LL_\Lambda := \LL_0 \cup \{\TT_\ga \mid \alpha <\Lambda\},$$ where $\TT_\ga$ are unary predicate symbols. The ordering $(\Lambda,<)$ determines a natural G\"odel numbering of all syntactic objects of $\LL_\Lambda$.

For each $\ga<\Lambda$ we interpret $\TT_\ga$ as the truth definition for the language $\LL_\ga$. So, we define an $\LL_{\alpha + 1}$-theory $\UTB_\alpha$ by the \emph{uniform Tarski biconditionals}:
\bi
\item[U1:] $\forall \vec x \left( \phi(\vec x) \eqv \TT_\ga(\gn{\phi(\num{\vec{x}})}) \right),$ for all $\LL_\ga$-formulas $\phi(\vec x)$ with all the free variables $\vec x$ shown;
\item[U2:] $\neg \TT(\num{n})$, for all $n$ such that $n$ is not a G\"odel number of an $\LL_\ga$-sentence.
\ei
Here, $\num{n}$ denotes the numeral for the natural number  $n$, and $\gn{\phi(\num{\vec x})}$ is an elementarily definable term representing the function mapping $n_1,\dots,n_k$ to the G\"odel number of $\phi(\num{n}_1,\dots,\num{n}_k)$, provided $\vec x=x_1,\dots,x_k$.

Also define
$$
\UTB_{<\alpha} := \bigcup_{\beta < \alpha} \UTB_\beta.
$$

\subsection{Hyperarithmetical hierarchy of formulas}

Let $\LL$ be a language extending $\LL_0$ by new predicate symbols.  $\Delta^{\LL}_0$ denotes the class of formulas obtained from atomic $\LL$-formulas by Boolean connectives and bounded quantifiers. The classes $\Pi_n^{\LL}$ and $\Sigma_n^{\LL}$ are defined from $\Delta^{\LL}_0$ in the usual way:  $\Pi_0^\LL=\Sigma_0^\LL=\Delta_0^\LL$, $\Pi_{n+1}^\LL=\{\al{\vec x}\phi(\vec x):\phi\in\gS^\LL_n\}$, and $\gS_{n+1}^\LL=\{\ex{\vec x}\phi(\vec x):\phi\in\Pi^\LL_n\}$.

Since one truth definition corresponds to $\gw$-many jumps, the ordinal notation system $(\Lambda,<)$ has to be extended to a slightly larger segment of ordinals up to $\gw(1+\gL)$, e.g., by encoding ordinals $\gw\ga+n$ as pairs $\la\ga,n\ra$. Then we introduce the following classes of formulas corresponding to the levels of the hyperarithmetical hierarchy up to $\gw(1+\gL)$ (the class $\Pi_{1+\ga}$ corresponds to $\Pi_1(\mathbf{0}^{(\ga)})$-sets):
\bi
\item $\Pi_n:=\Pi_n^{\LL_0}$ if $n<\gw$; \quad $\Pi_{\omega(1+\alpha) + n} := \Pi^{\LL_{\ga+1}}_{n+1}$ if $n<\gw$;
\item $\Pi_{< \alpha} := \bigcup_{\beta <\alpha} \Pi_\beta$. 
\ei

\subsection{Reflection principles}
Let $S$ be an r.e.\ extension of $\EA+\UTB_{<\Lambda}$ together with a fixed elementary formula defining the set of G\"odel numbers of its axioms in the standard model of arithmetic. Let $\Box_S$ denote the G\"odel's provability formula for $S$ (as defined, e.g., in \cite{Fef60}).

Suppose $\Gamma$ is a set of formulas in the language of $S$. By $\tRFN{\Gamma}{S}$ we denote the \emph{uniform reflection principle for $\Gamma$-formulas}, that is, the schema
$$\tRFN{\Gamma}{S}:\qquad \al{x}(\Box_S(\gn{\phi(\num{x})})\to \phi(x)), \quad \phi\in \Gamma.$$

More specifically, we define the reflection operators, for all $\ga<\gw(1+\gL)$, as follows:
\begin{eqnarray*}
\tR_{\alpha}(S) & := & \tRFN{\Pi_{1+\alpha}}{S}, \\
\tR_{<\ga}(S) & := & \text{$\tRFN{\Pi_{<\ga}}{S}$.}
\end{eqnarray*}

Theories axiomatized by transfinite iterations of reflection principles (along any elementary strict pre-wellordering $(\Omega,\prec)$) can be defined by formalizing the following equation in $S$:
\beq
\al{\gb\in\Omega}(\btR_\ga^\gb(S) \equiv S+ \textstyle{\bigcup}\{\tR_\ga(\btR_\ga^\gy(S)):\gy\prec \gb\}). \label{itr}
\eeq
Theories $\btR_\ga^\gb(S)$ satisfying \refeq{itr} are unique modulo provable equivalence in $S$ \cite{Bek18a}.
Similarly, one can define theories $\btR_{<\ga}^\gb(S)$ by transfinite iterations of the operators $\btR_{<\ga}:U\longmapsto S+\tR_{<\ga}(U)$.

Let $\EA^+$ denote $\EA$ together with the axiom asserting that the iterated exponentiation is total. For the main results of \cite{BekPakh}, and in this paper, we take $\IB:=\EA^++\UTB_{<\Lambda}$ as the base theory $S$. We will denote $\btR_{\ga}^\gb(\IB)$ and $\btR_{<\ga}^\gb(\IB)$ by $\btR_{\ga}^\gb$ and $\btR_{<\ga}^\gb$, respectively.

\section{Schmerl formulas}

One of the main results of \cite{BekPakh} is a conservation result relating mixed reflection principles to a transfinitely iterated reflection principle of complexity $\Pi_\ga$. Here we show that this result yields the relationships between the hierarchies of iterated reflection principles known as \emph{Schmerl formulas}. Such relationships first occurred in the work of Ulf Schmerl \cite{Schm,Sch82} and have been generalized in \cite{Bek99b}.

We consider a minor variant of the Veblen $\phi$ function:
$$
\bar\phi_\ga(\gb):=
\begin{cases} 0, \text{if $\gb=0$}, \\
\gw^\gb, \text{if $\ga=0$, $\gb\neq 0$}, \\
\phi_\ga(-1+\gb), \text{otherwise}.
\end{cases}
$$
Let $\Cr_\ga$ denote the range of $\bar\phi_\ga$, then $\Cr_{\ga+1}$ is the set of fixed points of $\bar\phi_\ga$ and, for limit ordinals $\gl$, $\Cr_\gl=\bigcup_{\ga<\gl}\Cr_\ga$. The functions $\bar\phi_\ga$ are increasing and continuous, and the sets $\Cr_\ga$ are closed and unbounded. In terms of the hierarchy of \emph{hyperation} functions introduced in \cite{FJ13a} we can express $\bar\phi_\ga(\gb)$ as $e^{\gw^\ga}(\gb)$.

We will use below specific ordinal notation systems $\Wo^\Lambda$ (and $\Wo^\Lambda_\ga$) associated with the reflection calculus $\Rcl$ (cf \cite[Section 6.2]{BekPakh}). These systems give notations to ordinals from an initial segment containing $\Lambda$ and on which the functions such as $+$ and $\phi_\ga$, for $\ga<\Lambda$, are well-defined.

Let $U\equiv_\ga V$ denote mutual conservativity of theories $U$ and $V$ w.r.t.\ $\Pi_{1+\ga}$-sentences.

\bt \label{Sch}
\benr
\item $\btR_{\ga+\gw^\gb}^\gy \equiv_\ga \btR_\ga^{\bar\phi_\gb(\gy)}$;
\item $\btR_{\ga+\gw^\gb}^{\gy_1}\cup \btR_\ga^{\gy_2+1} \equiv_\ga \btR_\ga^{\gy_2+1+\bar\phi_\gb(\gy_1)}$.
\eenr
\et

\bp\ (i) We assume $\Lambda$ so large that $\ga+\gw^\gb<\Lambda$ and $\gy$ belongs to the notation system $\Wo^\Lambda_{\ga+\gw^\gb}$ associated with the reflection calculus $\Rcl$. For each $\ga,\gb,\gy$ we can always select such a $\Lambda$ (and we may assume $\Lambda$ to be additively indecomposable).

This means that there is a word $C\in\Wo^{\Lambda}_{\ga+\gw^\gb}$ such that $o_{\ga+\gw^\gb}(C)=\gy$. By Theorem 8 of \cite{BekPakh} we have
$$C^*\equiv_{\ga+\gw^\gb} \btR_{\ga+\gw^\gb}^\gy.$$
Since $C\in \Wo^{\Lambda}_\ga$ the same result also yields
$$C^*\equiv_\ga \btR_{\ga}^{o_\ga(C)}.$$
Therefore, $\btR_{\ga+\gw^\gb}^\gy\equiv_\ga \btR_{\ga}^{o_\ga(C)}$ and it remains for us to compute $o_\ga(C)$.

Let $\nu\uparrow A$ denote the result of replacing each letter $x$ in $A$ by $\nu+x$. Similarly, $\nu\downarrow A$ denotes the result of replacing each letter $x$ in $A$ by $-\nu+x$. Due to the translation symmetry of $\Rcl$, which holds for additively indecomposable $\Lambda$, these maps provide obvious isomorphisms between the notation systems $(\Wo^{\Lambda}_\nu,<_\nu)$ and $(\Wo^{\Lambda}_0,<_0)$.

Let $D:= (\ga+\gw^\gb)\downarrow C=\gw^\gb\downarrow(\ga\downarrow C)$. We have $$o(D)=o_{\ga+\gw^\gb}(C)=\gy.$$ Then we obtain:
$$o_\ga(C)=o(\ga\downarrow C) = o(\gw^\gb\uparrow D) = \bar\phi_\gb(o(D))=\bar\phi_\gb(\gy),$$ by \cite[Lemma 17]{Bek05a} (see also \cite[Section 6.2]{BekPakh}).

(ii) Reasoning in a similar way consider $C_1\in\Wo^{\Lambda}_{\ga+\gw^\gb}$ and $C_2\in \Wo^{\Lambda}_{\ga}$ such that $o_{\ga+\gw^\gb}(C_1)=\gy_1$ and $o_\ga(C_2)=\gy_2$. Then
$$\btR_{\ga+\gw^\gb}^{\gy_1}\equiv_{\ga+\gw^\gb} C^*_1 \quad \text{ and } \quad \btR_{\ga}^{\gy_2}\equiv_\ga C^*_2.$$
The second equivalence yields $$\btR_\ga(\btR_{\ga}^{\gy_2})\equiv \btR_\ga(C^*_2).$$
Since $\tR_\ga$ has complexity $\Pi_{1+\ga}$ the first equivalence yields
$$\btR_{\ga+\gw^\gb}^{\gy_1}\cup \btR_\ga(\btR_{\ga}^{\gy_2})\equiv_{\ga+\gw^\gb} (C_1\land \ga C_2)^*.$$
We have that $C_1\land \ga C_2 =_\Rcl C_1\ga C_2$, since $C_1\in \Wo^{\Lambda}_{\ga+1}$. Moreover, $$o_\ga(C_1\ga C_2)=o_\ga(C_2)+1+\gw^{o_{\ga+1}(C_1)}=\gy_2+1+\bar\phi_\gb(\gy_1).$$
Then by \cite[Theorem 8]{BekPakh} we obtain
 $$(C_1\land \ga C_2)^*\equiv_\ga \btR_\ga^{\gy_2+1+\bar\phi_\gb(\gy_1)},$$
 as required.
\ep

We remark that formula (ii) can also be inferred from (i) by reflexive induction as in \cite[Lemma 7.3]{Bek18}.

\section{Conservativity spectra}
As before, we consider all ordinals to be represented in some notation system $\Wo^\Lambda$ for a suitably large ordinal $\Lambda$. We define the conservativity spectrum of a theory $S$ (of length $\gl$) as the $\gl$-sequence of ordinals $\ord_\ga(S)$, for all $\ga<\gl$, where
$$\ord_\ga(S):=\sup\{\gy\in \Wo^\Lambda: S\vdash \btR_\ga^\gy\}.$$
The conservativity spectrum is \emph{proper}, if the value $\ord_\ga(S)$ is in $\Wo^\Lambda$ for each $\ga<\gl$. We will tacitly assume all considered conservativity spectra to be proper.

Obviously, the sequence $\ord_\ga(S)$ is non-increasing with $\ga$.
The following theorem gives a stronger necessary condition for a $\gl$-sequence of ordinals to represent the conservativity spectrum of $S$ for some theory $S$.

\bt \label{spec} Suppose $f$ is the conservativity spectrum of $S$ of length $\gl$. For all $\ga,\gb$ such that $\ga+\gw^\gb<\gl$,
\benr
\item  $\ell(f(\ga))\geq f(\ga+1)$;
\item $\ell(f(\ga))\geq \bar\phi_\gb(f(\ga+\gw^\gb))$ if $\gb>0$.
\eenr
\et

\bp\
(i) Let $\gy:=f(\ga+1)$ and assume for a contradiction that $\gy>\ell f(\ga)$. In this case $f(\ga)\geq f(\ga+1)=\gy>0$ and we can write $f(\ga)=\ga_0+\gw^{\ell f(\ga)}$ for some $\ga_0$.

By the definition of spectrum $$S\vdash \btR_{\ga+1}^\gy\cup \btR_\ga^{\ga_0+1}.$$
By Theorem \ref{Sch} (ii) it follows that
$S\vdash  \btR_\ga^{\ga_0+1+\gw^\gy}$.
Hence, $\ord_\ga(S)\geq \ga_0+\gw^\gy >\ga_0+\gw^{\ell f(\ga)}=f(\ga)$, a contradiction.

(ii) Let $\gy:= f(\ga+\gw^\gb)$, for some $\gb>0$, and assume for a contradiction that $\ell f(\ga) < \bar\phi_\gb(\gy)$. Since $\gb>0$ we then also have $\bar\phi_\gb(\gy)= \gw^{\bar\phi_\gb(\gy)}>\gw^{\ell f(\ga)}$.

By the definition of spectrum
$$S\vdash \btR_{\ga+\gw^\gb}^\gy\cup \btR_\ga^{\ga_0+1}.$$
By Theorem \ref{Sch} (ii) it follows that
$S\vdash  \btR_\ga^{\ga_0+1+\bar\phi_\gb(\gy)}$.
Hence, $\ord_\ga(S)\geq \ga_0+\bar\phi_\gb(\gy) >\ga_0+\gw^{\ell f(\ga)}=f(\ga)$, a contradiction.
\ep

Fern\'andez and Joosten \cite[Proposition 5.2]{FJ13a} define the notion of \emph{$\ell$-sequence} as an ordinal sequence of length $\gl$ such that, for all $\xi<\zeta<\gl$,
$$\ell f(\xi)\geq \ell e^{\gw^{\ell\zeta}}f(\zeta).$$
As we recall, their function $e$ is such that $e^{\gw^\ga}(\gb)=\bar\phi_\ga(\gb)$.
Therefore, their condition is equivalent to demanding that $\ell f(\xi)\geq f(\zeta)$ if $\zeta=\xi+1$, and that
$$\ell f(\xi)\geq  \bar\phi_\gb(f(\zeta))$$
if $\zeta=\xi+\gw^\gb$ for some $\gb>0$. (If $\gb>0$ then $\bar\phi_\gb(f(\zeta))$ is a fixed point of $\ell$, hence the $\ell$ in front of $\bar\phi$ can be omitted.)
Hence, the necessary condition in Theorem \ref{spec} is equivalent to $f$ being an $\ell$-sequence.

\bcor The $\gl$-conservativity spectrum of any theory $S$ is an $\ell$-sequence.
\ecor

Before showing that every $\ell$-sequence is a conservativity spectrum of some theory let us notice a few properties of $\ell$-sequences.
Assume $f$ is an $\ell$-sequence and let $\gy_1:=f(\ga)$ and $\gy_2:=\min\{f(\gy):\gy<\ga\}$.
\bl \label{aless}
$\btR_\ga^{\gy_1}\cup \btR_{<\ga}^{\gy_2}\equiv_{<\ga} \btR_{<\ga}^{\gy_2}.$
\el

\bp\ If $\ga=\ga_0+1$, for some $\ga_0$, then $\gy_2=f(\ga_0)$, $\ell(\gy_2)\geq \gy_1$ and  we need to show
$$\btR_{\ga_0+1}^{\gy_1}\cup \btR_{\ga_0}^{\gy_2}\equiv_{\ga_0} \btR_{\ga_0}^{\gy_2}.$$
We may also assume $\gy_1>0$, otherwise the claim is trivial, hence $\gy_2\in\Lim$.

Consider any successor ordinal $\gd<\gy_2$. Then
$$\btR_{\ga_0+1}^{\gy_1}\cup \btR_{\ga_0}^{\gd}\equiv_{\ga_0} \btR_{\ga_0}^{\gd+\gw^{\gy_1}}.$$
Now we notice that $\gd+\gw^{\gy_1}\leq \gd +\gw^{\ell\gy_2}\leq \gy_2$. Hence, for any $\gd<\gy_2$,
$\btR_{\ga_0+1}^{\gy_1}\cup \btR_{\ga_0}^{\gd}$ is $\Pi_{1+\ga_0}$-conservative over
$\btR_{\ga_0}^{\gy_2},$ which proves the claim.

Assume $\ga=\ga_0+\gw^\gb$ with $\gb>0$. Pick any $\ga'<\ga$ such that $\ga'\geq \ga_0$ and $f(\ga')=\gy_2$. Then $\ga'+\gw^{\gb}=\ga$ and $\ell(\gy_2)\geq  \bar\phi_\gb(\gy_1)$.

Pick any successor ordinal $\gd<\gy_2$. We have:
$$\btR_{\ga}^{\gy_1}\cup \btR_{\ga'}^{\gd}\equiv_{\ga'} \btR_{\ga'}^{\gd+\bar\phi_{\gb}(\gy_1)}.$$
Since $\ell\gy_2\geq \bar\phi_\gb(\gy_1)$ we have $\gd+\bar\phi_\gb(\gy_1)\leq \gy_2$. Since this holds for all sufficiently large $\ga'<\ga$ and $\gd<\gy_2$, the claim follows.
\ep

\bt \label{spec-inv} Every $\ell$-sequence of length $\gl$ is a conservativity spectrum of some theory $S$.
\et
\bp\
Notice that every $\ell$-sequence $f$ is non-increasing and therefore has at most finitely many different values, say $\gy_1,\gy_2,\dots,\gy_n$ in the decreasing order. Let $$\ga_i=\min\{\ga: f(\ga)=\gy_{i+1}\}.$$ Then $f$ is constant $\gy_{i+1}$ on each interval $[\ga_i,\ga_{i+1})$, where $\ga_0=0$ and we put $\ga_{n}:=\gl$. Let
$$S_n:=\btR_{<\ga_{n}}^{\gy_n} \cup \btR_{<\ga_{n-1}}^{\gy_{n-1}}\cup \cdots \cup \btR_{<\ga_{1}}^{\gy_1}.$$
We claim that the conservativity spectrum of $S_n$ coincides  with $f$. To show this we need two more lemmas.

\bl \label{const} Suppose $\ga_{i-1}\leq\ga<\ga_i$. Then,
$\btR_{\ga_{i}}^{\gy_{i}}\equiv_{\ga} \btR_{\ga}^{\gy_{i}}.$
\el
\bp\ We can write $\ga_i:=\ga+\gw^{\gb_1}+\cdots +\gw^{\gb_k}$. Let $\bar\ga_j:= \ga+\gw^{\gb_1}+\cdots +\gw^{\gb_j}$ with $\bar\ga_0:=\ga$. Since $f$ is an $\ell$-sequence and $f(\bar\ga_j)=\gy_i$, we have $\gy_i\geq \ell(\gy_i)\geq \bar\phi_{\gb_j}(\gy_i)$, hence $\gy_i$ is a fixed point of $\phi_{\gb_j}$, for each $j$. Then by induction on $j=k,\dots,0$ from Theorem \ref{Sch} (i) we obtain that
$$\btR_{\ga_{i}}^{\gy_{i}}\equiv_{\bar\ga_{j}} \btR_{\bar\ga_{j}}^{\gy_{i}}.$$ The claim follows from this for $j=0$. \ep

\bl \label{cons} For each $i< n$, $S_{i+1}\equiv_{<\ga_{i}} S_i$.
\el
\bp\ Firstly, by Lemma \ref{const},  $$\btR_{<\ga_{i+1}}^{\gy_{i+1}}\equiv_{<\ga_{i+1}}\btR_{\ga_{i+1}}^{\gy_{i+1}}\equiv_{\ga_i} \btR_{\ga_{i}}^{\gy_{i+1}}.$$
Since $S_{i}$ is a set of formulas of complexity $\Pi_{<\ga_i}$, it follows that
$$S_{i+1}\equiv S_i\cup \btR_{<\ga_{i+1}}^{\gy_{i+1}} \equiv_{\ga_i} S_i\cup \btR_{\ga_{i}}^{\gy_{i+1}}
 \equiv  S_{i-1}\cup \btR_{<\ga_{i}}^{\gy_{i}}\cup \btR_{\ga_{i}}^{\gy_{i+1}}.$$
Since $S_{i-1}$ has complexity $\Pi_{<\ga_{i-1}}$, Lemma \ref{aless} implies that
$$S_{i-1}\cup \btR_{<\ga_{i}}^{\gy_{i}}\cup \btR_{\ga_{i}}^{\gy_{i+1}} \equiv_{<\ga_i}  S_{i-1}\cup \btR_{<\ga_{i}}^{\gy_{i}} \equiv S_i.$$
\ep

Now we prove by induction on $i=0,\dots, n$ that the $\ga_i$-conservativity spectrum of $S_i$ coincides with $f\rst \ga_i$. Theorem \ref{spec-inv} is this statement for $i=n$.
For $i=0$ the statement holds trivially.

We assume it holds for $i$ and prove it for $i+1$. By Lemma \ref{cons}, $\ord_\ga(S_{i+1})=\ord_\ga(S_i)$, for all $\ga<\ga_i$.
So we consider an $\ga$ such that $\ga_i\leq\ga<\ga_{i+1}$. Then, since the complexity of $S_i$ is $\Pi_{<\ga_i}$ and $S_i$ is sound, $$\ord_\ga(S_{i+1})=\ord_\ga(S_{i}\cup \btR_{<\ga_{i+1}}^{\gy_{i+1}})=
\ord_\ga(\btR_{<\ga_{i+1}}^{\gy_{i+1}}).$$
Finally, by Lemma \ref{const}, $$\btR_{<\ga_{i+1}}^{\gy_{i+1}}\equiv_\ga \btR_{\ga}^{\gy_{i+1}}.$$
Hence, $\ord_\ga(\btR_{<\ga_{i+1}}^{\gy_{i+1}})=\gy_{i+1}=f(\ga)$. Theorem \ref{spec-inv} is proved.
\ep

\bibliographystyle{plain}
\bibliography{ref-all2}

\end{document}